\documentclass[11.05pt,a4paper]{amsart}
\usepackage{amssymb,amsmath}
\usepackage{color}
\usepackage{latexsym}
\usepackage{amsthm,amsfonts,amssymb,mathrsfs}
\usepackage{rotating}
\usepackage[leqno]{amsmath}
\usepackage{xspace}
\usepackage[all]{xy}
\usepackage{longtable}

\headsep=1cm \numberwithin{equation}{section}
\newtheorem{theorem}{Theorem}[section]
\newtheorem{lemma}[theorem]{Lemma}
\newtheorem{proposition}[theorem]{Proposition}
\newtheorem{corollary}[theorem]{Corollary}

\theoremstyle{definition}
\newtheorem{definition}[theorem]{Definition}
\theoremstyle{remark}
\newtheorem{remark}[theorem]{Remark}

\newcommand{\Ass}{\operatorname{Ass}}

\newcommand{\sdepth}{\operatorname{sdepth}}

\newcommand{\Mon}{\operatorname{Mon}}

\newcommand{\G}{\operatorname{G}}

\newcommand{\Ht}{\operatorname{ht}}

\newcommand{\reg}{\operatorname{reg}}

\newcommand{\Proj}{\operatorname{Proj}}

\newcommand{\depth}{\operatorname{depth}}

\newcommand{\fp}{\frak{p}}

\newenvironment{prf}[1][Proof]{\begin{proof}[\bf #1]}{\end{proof}}
\begin{document}

\author[S. Bandari, K. Divaani-Aazar and A. Soleyman Jahan]{Somayeh Bandari,
Kamran Divaani-Aazar and Ali Soleyman Jahan}
\title[Almost complete intersections and ...]
{Almost complete intersections and Stanley's conjecture}

\address{S. Bandari, Department of Mathematics, Alzahra
University, Vanak, Post Code 19834, Tehran, Iran.}
\email{somayeh.bandari@yahoo.com}

\address{K. Divaani-Aazar, Department of Mathematics, Alzahra
University, Vanak, Post Code 19834, Tehran, Iran-and-School of
Mathematics, Institute for Research in Fundamental Sciences (IPM),
P.O. Box 19395-5746, Tehran, Iran.}
\email{kdivaani@ipm.ir}

\address{A. Soleyman Jahan, Department of Mathematics, University of Kurdistan,
Post Code 66177-15175, Sanandaj, Iran.}
\email{solymanjahan@gmail.com}

\subjclass[2010]{13F20; 05E40; 13F55.}

\keywords {Almost complete intersection monomial ideals; clean;  locally complete intersection monomial
ideals; pretty clean.\\
The research of the second and third authors are supported by grants from
IPM (no. 90130212 and no. 90130062, respectively).}

\begin{abstract} Let $K$ be a field and $I$ a monomial ideal of the polynomial ring $S=K[x_1,\ldots, x_n]$.
We show that if either: 1) $I$ is almost complete intersection, 2) $I$ can be generated by less than four
monomials; or 3) $I$ is the Stanley-Reisner ideal of a locally complete intersection simplicial complex on
$[n]$, then Stanley's conjecture holds for $S/I$.
\end{abstract}

\maketitle

\section{Introduction}

Throughout this paper, let $K$ be a field and $I$ a monomial ideal of the polynomial ring
$S=K[x_1,\ldots, x_n]$.

A decomposition of $S/I$ as direct sum of $K$-vector spaces of the form $\mathcal D: S/I=
\bigoplus_{i=1}^ru_iK[Z_i]$, where $u_i$ is a monomial in $S$ and $Z_i\subseteq \{x_1,\ldots, x_n\}$,
is called a {\em Stanley decomposition} of $S/I$.
The number $\sdepth \mathcal D:=\min\{|Z_i|: i=1,\ldots, r\}$ is called {\em Stanley depth}
of $\mathcal D$. The {\em Stanley depth} of $S/I$ is defined to be $$\sdepth S/I:=\max \{\sdepth
\mathcal D: \mathcal D \ \text{is a Stanley decomposition of} \  S/I\}.$$
Stanley conjectured \cite{St} that $\depth S/I\leq \sdepth S/I$. This conjecture is known
as Stanley's conjecture. Recently, this conjecture was extensively examined by several
authors; see e.g. \cite{A1}, \cite{A2}, \cite{HP},
\cite{HSY}, \cite{P}, \cite{R}, \cite{S2} and \cite{S3}. On the other hand, the present
third author \cite{S2} conjectured that there always exists a Stanley decomposition
$\mathcal D$ of $S/I$ such that the degree of each $u_i$ is at most $\reg S/I$. We refer
to this conjecture as $h$-regularity conjecture. It is known that for square-free
monomial ideals, these two conjectures are equivalent.  Our main aim in this paper is
to determine some classes of monomial ideals such that these conjectures are true for them.

A basic fact in commutative algebra says that there exists a finite chain $$\mathcal{F}:
I=I_0\subset I_1\subset \cdots \subset I_r=S$$ of monomial ideals such that $I_i/I_{i-1}
\cong S/\frak p_i$ for monomial prime ideals $\frak p_i$ of $S$. Dress \cite{D} called the
ring $S/I$ {\em clean} if there exists a chain $\mathcal{F}$ such that all the $\frak p_i$
are minimal prime ideals of $I$. By \cite[Proposition 2.2]{HSY} if $I$ is complete intersection,
then the ring $S/I$ is clean. Lemmas \ref{almost} and 2.8 provide two other classes of clean rings.

Herzog and Popescu \cite{HP} called the ring $S/I$ {\em pretty clean} if there exists a chain
$\mathcal{F}$ such that for all $i<j$ for which $\frak p_i\subseteq \frak p_j$, it follows that
$\frak p_i=\frak p_j$. Obviously, cleanness implies pretty cleanness and when $I$ is square-free,
it is known that these two concepts coincide; see \cite[Corollary 3.5]{HP}.

If $S/I$ is pretty clean, then $S/I$ is sequentially Cohen-Macaulay and depth of $S/I$ is equal
to the minimum of the dimension of $S/\frak p$, where $\frak p\in\Ass_SS/I;$ see \cite{S1} for
an easy proof. If $S/I$ is pretty clean, then \cite[Theorem 6.5]{HP} asserts that Stanley's
conjecture holds for $S/I$. In fact, if $S/I$ is pretty clean, then \cite[Proposition 1.3]{HVZ}
yields that $\depth S/I=\sdepth S/I$. Also if $S/I$ is pretty clean, then by \cite[Theorem 4.7]{S2}
$h$-regularity conjecture holds for $S/I$.

We prove that if the monomial ideal $I$ is either almost complete intersection or it can be generated
by less than four monomials, then $S/I$ is pretty clean. Thus, for such monomial ideals both Stanley's
and $h$-regularity conjectures hold. Also, we show that if $I$ is the Stanley-Reisner ideal of
a locally complete intersection simplicial complex on $[n]$, then $S/I$ satisfies Stanley's conjecture.

\section{Main Results}

A {\em simplicial complex} $\Delta$ on $[n]:=\{1,\ldots, n\}$ is a collection of subsets of
$[n]$ with the property that if $F\in \Delta$, then all subsets of $F$ are also in $\Delta$. Any
singleton element of $\Delta$ is called a {\em vertex}. An element of $\Delta$ is called a
{\em face} of $\Delta$ and the maximal faces of $\Delta$, under inclusion, are called
{\em facets}. We denote by $\mathcal F(\Delta)$ the set of all facets
of $\Delta$. The {\em dimension} of a face $F$ is defined as $\dim F=|F|-1$, where $|F|$ is
the number of elements of $F$. The dimension of the simplicial complex $\Delta$ is the maximal
dimension of its facets. A simplicial complex $\Delta$ is called {\em pure} if all facets of
$\Delta$ have the same dimension. We denote the simplicial complex $\Delta$ with facets
$F_1,\ldots, F_t$ by $\Delta=\langle F_1,\ldots, F_t\rangle$. According to Bj\"{o}rner and Wachs
\cite{BW},  a simplicial complex $\Delta$ is said to be {\em (non-pure) shellable} if there
exists an order $F_1,\ldots, F_t$ of the facets of $\Delta$ such that
for each $2\leq i\leq t$, $\langle F_1,\ldots, F_{i-1}\rangle \cap
\langle F_i \rangle$ is a pure $(\dim F_i-1)$-dimensional simplicial complex. If $\Delta$ is a
simplicial complex on $[n]$, then the {\em Stanley-Reisner ideal} of $\Delta$, $I_\Delta$, is the
square-free monomial ideal generated by all monomials $x_{i_1}x_{i_2}\ldots x_{i_t}$ such
that $\{i_1,i_2,\ldots, i_t\} \not\in \Delta$. The {\em Stanley-Reisner ring} of $\Delta$ over
the field $K$ is the $K$-algebra $K[\Delta]:=S/I_\Delta$. Any square-free monomial ideal $I$ is the
Stanley-Reisner ideal of some simplicial complex $\Delta$ on $[n]$. If $\mathcal F(\Delta)=
\{F_1,\ldots, F_t\}$, then $I_\Delta=\bigcap_{i=1}^t \frak p_{F_i}$, where $\frak p_{F_i}:=
(x_j:j\not\in F_i)$; see \cite[Theorem 5.1.4]{BH}.

Recall that the {\em Alexander dual} $\Delta^\vee$  of a simplicial complex $\Delta$ is
the simplicial complex whose faces are $\{[n]\backslash F|F\notin \Delta\}$. Let $I$ be
a square-free monomial ideal of $S$. We denote by $I^\vee,$ the square-free monomial ideal
which is generated by all monomials $x_{i_1}\cdots x_{i_k}$, where $(x_{i_1},\ldots,
x_{i_k})$ is a minimal prime ideal of $I$. It is easy to see that for any simplicial complex
$\Delta$, one has $I_{\Delta^\vee}=(I_{\Delta})^\vee$. A monomial ideal $I$ of $S$ is said to
have {\em linear quotients} if there exists an order $u_1,\ldots, u_m$ of $\G(I)$ such that
for any $2\leq i\leq m$, the ideal $(u_1,\ldots, u_{i-1}):_Su_i$ is generated by a subset of the
variables.

\begin{lemma}\label{vech} Let $I$ be a square-free monomial ideal of $S$. Then $S/I$ is clean
if and only if $I^\vee$ has linear quotients.
\end{lemma}

\begin{prf} Dress \cite[Theorem on page 53]{D} proved that a simplicial complex $\Delta$
is (non-pure) shellable if and only if $K[\Delta]$ is a clean ring. On the other hand, by
\cite[Theorem 1.4]{HHZ}, a simplicial complex $\Delta$ is (non-pure) shellable if and only
if $I_{\Delta^\vee}$ has linear quotients. Combining these facts, yields our claim.
\end{prf}

\begin{lemma}\label{Dr} Let $I$ and $J$ be two monomial ideals of $S$. Assume that $I=uJ$ for
some monomial $u$ in $S$ and $\Ht J\geq 2$. If $S/J$ is pretty clean, then $S/I$ is pretty
clean too.
\end{lemma}

\begin{prf} With the proof of \cite[Lemma 1.9]{S3}, the claim is immediate.
\end{prf}

In what follows for a monomial ideal $I$ of $S$, we denote the number of elements
of $\G(I)$ by $\mu(I)$.

\begin{definition} A monomial ideal $I$ of $S$ is said to be {\em almost complete
intersection} if $\mu(I)=\Ht I+1$.
\end{definition}

\begin{lemma}\label{almost} Let $I$ be an almost complete intersection square-free
monomial ideal of $S$. Then $S/I$ is clean.
\end{lemma}

\begin{prf} The claim is obvious when $\Ht I=0$. Let $\Ht I=1$. Then $I=(u_1,u_2)$ for
some monomials $u_1$ and $u_2$. We can write $I$ as $I=u(u'_1,u'_2)$, where $u=\gcd(u_1,u_2)$
and $u'_1,u'_2$ are monomials forming a regular sequence on $S$. So in this case, the claim
is immediate by Lemma \ref{Dr} and \cite[Proposition 2.2]{HSY}. Now, assume that $h:=\Ht I\geq 2$.
By \cite[Theorem 4.4]{KTY} $I$ can be written in one of the following forms, where
$A_1,A_2,\ldots, B_1,B_2, \ldots$ are non-trivial square-free monomials which are pairwise
relatively prime, and $p,p'$ are integers with $2\leq p\leq h$ and $1\leq p'\leq h$.
\begin{enumerate}
\item[1)] $I_1=(A_1B_1,A_2B_2, \ldots , A_pB_p,A_{p+1}, \ldots, A_h, B_1B_ 2\cdots B_p).$
\item[2)] $I_2=(A_1B_1,A_2B_2, \ldots , A_{p'}B_{p'},A_{p'+1}, \ldots, A_h,A_{h+1}B_1B_2
\cdots B_{p'})$.
\item[3)] $I_3=(B_1B_2,B_1B_3,B_2B_3,A_4, \ldots ,A_{h+1})$.
\item[4)] $I_4=(A_1B_1B_2,B_1B_3,B_2B_3,A_4, \ldots ,A_{h+1})$.
\item[5)] $I_5=(A_1B_1B_2,A_2B_1B_3,B_2B_3,A_4, \ldots ,A_{h+1})$.
\item[6)] $I_6=(A_1B_1B_2,A_2B_1B_3,A_3B_2B_3,A_4, \ldots ,A_{h+1})$.
\end{enumerate}
Let $I=I_1$. Since $A_1,A_2, \ldots , A_p,A_{p+1}, \ldots, A_h,B_1,B_2, \ldots ,B_p$
are pairwise relatively prime, it turns out that $A_{p+1}, \ldots, A_h$ is a regular sequence on
$S/(A_1B_1,A_2B_2, \ldots , A_pB_p, B_1B_ 2\cdots B_p)$. So, in view of
\cite[Theorem 2.1]{R}, we may and do assume that $I=(A_1B_1,A_2B_2, \ldots , A_pB_p, B_1B_ 2
\cdots B_p)$. Next, we are going to show that $I$ is of forest type. Let $G$ be a subset of
$\{A_1B_1,A_2B_2, \ldots , A_pB_p, B_1B_ 2\cdots B_p\}$ with at least two elements. If
$B_1B_ 2\cdots B_p\notin G$, then any $a\in G$ can be taken as a leaf and any $b\in G$ different
from $a$ can be taken as a branch for this leaf. If $B_1B_ 2\cdots B_p\in G$, then any
$a\in G$ different from $B_1B_ 2\cdots B_p$ can be taken as a leaf and then $B_1B_ 2\cdots B_p$
is a branch for this leaf. So, $I$ is of forest type. Thus, since $I$ is square-free, by
\cite[Theorem 1.5]{SZ}, we obtain that $S/I$ is clean. By the similar argument, one can see
that if $I=I_2$, then  $S/I$ is clean. Set $$J:=(C_1B_1B_2,C_2B_1B_3,
C_3B_2B_3,A_4, \ldots, A_{h+1}),$$ where $C_i$ is either $A_i$ or 1 for each $i=1,2,3$.
Since each of $I_3$, $I_4$, $I_5$ and $I_6$ are the particular cases of the ideal $J$, we
can finish the proof by showing that $S/J$ is clean.
Since, by the assumption $A_4,\ldots, A_{h+1},B_1,B_2,B_3,C_1,C_2,C_3$ are pairwise relatively
prime, it follows that $A_4,\ldots, A_{h+1}$ is a regular sequence on $S/(C_1B_1B_2,C_2B_1B_3,
C_3B_2B_3)$. So by
\cite[Theorem 2.1]{R}, we can assume that $J=(C_1B_1B_2,C_2B_1B_3,C_3B_2B_3)$. Set $T:=k[u,v,w,x,y,z]$
and $L:=(uxy,vxz,wyz)$. Since $B_1,B_2,B_3,C_1,C_2,C_3$ is  a regular sequence on $S$, by
\cite[Proposition 3.3]{HSY}, the cleanness of $T/L$ implies the cleanness of $S/J$. So, by Lemma
\ref{vech}, it is enough to prove that $L^\vee$ has linear quotients. As $$L=(x,y)\cap (x,z)\cap (x,w)
\cap (y,z)\cap (y,v) \cap(z,u) \cap(u,v,w),$$ one has $L^\vee=(xy,xz,xw,yz,yv,zu,uvw),$ which clearly
has linear quotients by the given order.
\end{prf}

Let $u=\prod_{i=1}^nx_i^{a_i}$ be a monomial in $S=K[x_1,\ldots, x_n]$. Then
$$u^p:=\textstyle\prod\limits_{i=1}^n\textstyle\prod\limits_{j=1}^{a_i}x_{i,j}
\in K[x_{1,1},\ldots, x_{1,a_1},\ldots, x_{n,1},\ldots, x_{n,a_n}]$$ is called the
{\em polarization} of $u$. Let  $I$ be a monomial ideal of $S$ with $\G(I)=
\{u_1,\ldots, u_m\}$. Then the ideal $I^p:=(u_1^p,\ldots, u^p_m)$ of $T:=K[x_{i,j}:i=1,
\ldots, n,  j=1,\ldots ,a_i]$ is called the {\em polarization} of $I$.
\cite[Theorem 3.10]{S3} implies that $S/I$ is
pretty clean if and only if $T/I^p$ is clean.

Recently, Cimpoea\c{s} \cite{C1} proved that if  $I$ is an almost complete
intersection monomial ideal of $S$, then Stanley's conjecture holds for $S/I$.
The next result shows that in this case $S/I$ is even pretty clean.

\begin{theorem}\label{almost p} Let $I$ be an almost complete intersection monomial
ideal of $S$. Then $S/I$ is pretty clean.
\end{theorem}

\begin{prf} From \cite[Proposition 2.3]{F}, one has $\Ht I=\Ht I^p$. On the
other hand $\mu(I)=\mu(I^p)$, and so $I^p$ is an almost complete intersection
square-free monomial ideal of $T$. Hence, by Lemma \ref{almost}, the ring
$T/I^p$ is clean. Now, \cite[Theorem 3.10]{S3} implies that $S/I$ is pretty clean,
as desired.
\end{prf}

In the situation of Theorem \ref{almost p}, there is no need that $S/I$  is clean.
For instance, although $(x^2,xy)$ is an almost complete intersection monomial ideal,
the ring $k[x,y]/(x^2,xy)$ is not clean.

In \cite[Theorem 2.3]{C2}, it is shown that if $I$ is a monomial ideal of $S$ with
$\mu(I)\leq 3$, then Stanley's  conjecture holds for $S/I$. The next result
extends this fact.

\begin{corollary} Let $I$ be a monomial ideal of $S$. If $\mu(I)\leq 3$, then $S/I$
is pretty clean.
\end{corollary}

\begin{prf} Clearly, we may assume that $I$ is non zero. Assume that $\mu(I)=3$ and
$\Ht I=1$. Then $I=uJ$, where $u$ is a monomial in $S$ and $J$ is a monomial ideal of
$S$ with $\mu(J)=3$ and $\Ht J \geq 2$. By Lemma \ref{Dr}, it is
enough to prove that $S/J$ is pretty
clean. If $\Ht J=2$, then $\mu(J)=\Ht J+1$, and so by Theorem \ref{almost p}, $S/J$ is
pretty clean. If $\Ht J=3$, then $J$ is complete intersection, and hence by
\cite[Proposition 2.2]{HSY}, $S/J$ is pretty clean.

Since $0<\Ht I\leq \mu(I)$, in all other cases, it follows that $I$ is either complete
intersection or almost complete intersection. Thus, the proof is completed by
\cite[Proposition 2.2]{HSY} and Theorem \ref{almost p}.
\end{prf}

\begin{definition} (\cite[Definition 1.1 and Lemma 1.2]{TY}) A simplicial complex $\Delta$
on $[n]$ is said to be {\em locally complete intersection} if $\{\{1\},\{2\}, \ldots, \{n\}\}
\subseteq  \Delta$ and $(I_{\Delta})_{\fp}$ is a complete intersection ideal of $S_{\fp}$ for
all $\frak p\in \Proj S/I$.
\end{definition}

A simplicial complex $\Delta$ is said to be {\em connected} if for any two facets
$F$ and $G$ of $\Delta$, there exists a sequence of facets $F=F_0, F_1, \ldots,
F_{q-1}, F_q=G$ such that $F_i\cap F_{i+1}\neq\emptyset$ for all $0\leq i<q$.
Also, a simplicial complex $\Delta$ on $[n]$ is said to be $n$-{\em pointed path}
(resp. $n$-{\em gon}) if $n\geq 2$ (resp. $n\geq 3$) and, after a suitable change
of variables, $$\mathcal{F}(\Delta)=\{\{i,i+1\}|1\leq i<n\}$$ (resp.
$$\mathcal{F}(\Delta)=\{\{i,i+1\}|1\leq i<n\}\cup \{\{n,1\}\}).$$ Clearly, any
$n$-pointed path (resp. $n$-gon) is one-dimensional and pure.

Let $\Delta$ be a connected simplicial complex on $[n]$ which is locally complete intersection.
Then, it is known that $\Delta$ is shellable; see e.g. \cite[Proposition 1.11 and Theorem 1.5]{TY}.
Hence, by \cite[Theorem on page 53]{D} it turns out that $S/I_{\Delta}$ is clean. So, we record
the following:

\begin{lemma} Let $\Delta$ be a connected simplicial complex on $[n]$ which is locally complete
intersection. Then $S/I_{\Delta}$ is clean.
\end{lemma}

Let $\Delta$ be as in Lemma 2.8. Then $S/I_{\Delta}$ is clean, and so \cite[Theorem 6.5]{HP} implies
that $S/I_{\Delta}$ satisfies Stanley's conjecture. In Theorem \ref {second}, we prove that the later
assertion holds without assuming that $\Delta$ is connected.

\begin{proposition}\label{main} Let $I\subset S_1=K[x_1,\ldots,x_m]$,
$J\subset S_2= K[x_{m+1}, \ldots, x_n]$ be two monomial ideals and
$S=K[x_1,\ldots, x_m,x_{m+1},\ldots, x_n]$. Assume that $\depth S_1/I>0$ and
$\depth S_2/J>0$. Then Stanley's  conjecture holds for
$S/(I,J,\{x_ix_j\}_{1\leq i\leq m,m+1\leq j\leq n})$.
\end{proposition}

\begin{prf} For convenience, we set $Q_1:=(x_1,\ldots, x_m)$, $Q_2:=(x_{m+1},\ldots, x_n)$
and $Q:=(x_ix_j)_{1\leq i\leq m, m+1\leq j\leq n}$. So, $Q=Q_1\cap Q_2$. Since
$I\subseteq Q_1$ and $J\subseteq Q_2$, it follows that $$(I,J,Q)=(I,J,Q_1)\cap (I,J,Q_2)=
(J,Q_1)\cap (I,Q_2).$$ By the assumption, we have $Q_1\notin \Ass_{S_1}S_1/I$ and $Q_2\notin
\Ass_{S_2}S_2/J$. Hence
$$(x_1,\ldots, x_m,x_{m+1},\ldots, x_n)\notin \Ass_SS/(I,Q_2)$$ and
$$(x_1,\ldots, x_m,x_{m+1},\ldots, x_n)\notin \Ass_SS/(J,Q_1),$$
and so $$\depth(\frac{S}{(J,Q_1)}\oplus\frac{S}{(I,Q_2)})>0=\depth(\frac{S}{Q_1+Q_2}).$$
Now, in view of the exact sequence $$0\rightarrow \frac{S}{(J,Q_1)\cap(I,Q_2)}
\rightarrow \frac{S}{(J,Q_1)}\oplus\frac{S}{(I,Q_2)}\rightarrow\frac{S}{Q_1+Q_2}
\rightarrow 0,$$  \cite[Lemma 1.3.9]{V} implies that $$\depth(\frac{S}{(I,J,Q)})=
\depth(\frac{S}{(J,Q_1)\cap(I,Q_2)})=1.$$
Now the proof is complete, because \cite[Theorem 2.1]{C2} yields that for any monomial
ideals $L$ of $S$ if $\depth S/L\leq 1$, then Stanley's  conjecture holds for $S/L$.
\end{prf}

\begin{corollary}\label{first} Let $\Delta_1$ and $\Delta_2$ be two non-empty disjoint
simplicial complexes and $\Delta:=\Delta_1\cup \Delta_2$. Then Stanley's  conjecture holds
for $S/I_{\Delta}$.
\end{corollary}

\begin{prf} For two natural integers $m<n$, we may assume that $\Delta_1$ and $\Delta_2$
are simplicial complexes on $[m]$ and $\{m+1,\ldots,n\}$,  respectively. Then
$K[\Delta_1]=K[x_1,\ldots, x_m]/I_{\Delta_1}$ and
$K[\Delta_2]=K[x_{m+1},\ldots, x_n]/I_{\Delta_2}$, and so $$K[\Delta]=K[x_1,\ldots, x_m,x_{m+1},
\ldots, x_n]/(I_{\Delta_1},I_{\Delta_2}, \{x_ix_j\}_{1\leq i\leq m, m+1\leq j\leq n}).$$
We claim that $\depth(K[x_1,\ldots, x_m]/I_{\Delta_1})>0$ and
$\depth(K[x_{m+1},\ldots, x_n]/I_{\Delta_2})>0$.
Because if for example  $\depth(K[x_1,\ldots, x_m]/I_{\Delta_1})=0$, then
$I_{\Delta_1}=(x_1,\ldots, x_m)$. But, this implies that $\Delta_1=\emptyset$ which contradicts
our assumption on $\Delta_1$.  Now, the claim is immediate by Proposition \ref{main}.
\end{prf}

\begin{theorem}\label{second}  Let $\Delta$ be a locally complete intersection simplicial
complex on $[n]$. Then Stanley's conjecture holds for $S/I_{\Delta}$.
\end{theorem}

\begin{prf} If $\Delta$ is connected, then Lemma 2.8 yields the claim. Otherwise,
by \cite[Theorem 1.15]{TY}, $\Delta$ is the disjoint union of finitely many non-empty
simplicial complexes. So, in this case the assertion follows by Corollary \ref{first}.
\end{prf}

In \cite[Corollary 4.3]{HP} it is shown that if $S/I$ is pretty clean, then it is sequentially
Cohen-Macaulay. In \cite{S1} this fact is reproved by a different argument and, in addition, it
is shown that depth of $S/I$ is equal to the minimum of the dimension of $S/\frak p$, where
$\frak p\in\Ass_SS/I$.  Also if $S/I$ is pretty clean, then by \cite[Theorem 4.7]{S2}
$h$-regularity conjecture holds for $S/I$. This implies part a) of the following remark.

\begin{remark} Let $I$ be a monomial ideal of $S$.
\begin{enumerate}
\item[a)] Assume that either:\\
i) $I$ is almost complete intersection,\\
ii) $\mu(I)\leq 3$; or\\
iii) $I$ is the Stanley-Reisner ideal of a connected simplicial complex on $[n]$ which is
locally complete intersection.\\ Then both Stanley's  and
$h$-regularity conjectures hold for $S/I$. Also, in each of these cases $S/I$ is sequentially
Cohen-Macaulay and $\depth S/I=\min\{\dim S/\frak p|\frak p\in \Ass_SS/I\}$.
\item[b)] We know that if $S/I$ is pretty clean, then Stanley's  conjecture holds for $S/I$.
By using Corollary \ref{first}, we can provide an example of a monomial ideal $I$ of $S$ such that
Stanley's  conjecture holds for $S/I$, while it is not pretty clean. To this end, let
$\Delta_1$, $\Delta_2$ and $\Delta$ be as in Corollary \ref{first} and $\dim \Delta_i>0$, $i=1,2$.
Evidently, $\Delta$ is not shellable, and so \cite[Theorem on page 53]{D} implies that
$S/I_{\Delta}$ is not pretty clean. On the other hand, Stanley's  conjecture holds for
$S/I_{\Delta}$ by Corollary \ref{first}.
\end{enumerate}
\end{remark}



\begin{thebibliography}{99}

\bibitem[A1]{A1}{J. Apel}, {\it On a conjecture of R. P. Stanley. II. Quotients modulo
monomial ideals}, J. Algebraic Combin., {\bf 17}(1), (2003), 57-74.

\bibitem[A2]{A2}{J. Apel}, {\it On a conjecture of R. P. Stanley. I. Monomial ideals},
J. Algebraic Combin., {\bf 17}(1), (2003), 39-56.

\bibitem[BW]{BW}{A. Bj\"{o}rner and M. Wachs}, {\it Shellable nonpure complexes and posets. I},
Trans. Amer. Math. Soc., {\bf 348}(4), (1996), 1299-1327.

\bibitem[BH]{BH}{W. Bruns and J. Herzog}, {\it Cohen Macaulay rings},  Cambridge Studies in
Advanced Mathematics, {\bf 39}, Cambridge University Press, Cambridge, 1993.

\bibitem[C1]{C1}{M. Cimpoea\c{s}}, {\it The Stanley conjecture on monomial almost complete
intersection ideals}, Bull. Math. Soc. Sci. Math. Roumanie (N.S.), {\bf 55(103)}(1), (2012), 35-39.

\bibitem[C2]{C2}{M. Cimpoea\c{s}}, {\it Stanley depth of monomial ideals with small number of
generators}, Cent. Eur. J. Math., {\bf 7}(3), (2009), 629-634.

\bibitem[D]{D}{A. Dress}, {\it A new algebraic criterion for shellability},
Beitr$\ddot{a}$ge Algebra Geom., {\bf 34}(1), (1993), 45-55.

\bibitem[F]{F}{S. Faridi}, {\it Monomial ideals via square-free monomial ideals},
Commutative algebra, 85-114, Lect. Notes Pure Appl. Math., {\bf 244}, Chapman \& Hall/CRC, Boca
Raton, FL, (2006).

\bibitem[HHZ]{HHZ}{J. Herzog, T. Hibi and X. Zheng}, {\it Dirac's theorem on chordal graphs
and Alexander duality}, European J. Combin., {\bf 25}(7), (2004), 949-960.

\bibitem[HP]{HP}{J. Herzog and D. Popescu}, {\it Finite filtrations of modules and shellable
multicomplexes}, Manuscripta Math., {\bf 121}(3), (2006), 385-410.

\bibitem[HSY]{HSY}{J. Herzog, A. Soleyman Jahan and S. Yassemi}, {\it Stanley decompositions and
partitionable simplicial complexes}, J. Algebraic Combin., {\bf 27}(1), (2008), 113-125.

\bibitem[HVZ]{HVZ}{J. Herzog, M. Vladoiu and X. Zheng}, {\it How to compute the Stanley depth
of a monomial ideal}, J. Algebra, {\bf 322}(9), (2009), 3151-3169.

\bibitem[KTY]{KTY}{K. Kimura, N. Terai and K. Yoshida}, {\it Arithmetical rank of square-free
monomial ideals of small arithmetic degree},  J. Algebraic Combin., {\bf 29}(3), (2009), 389-404.

\bibitem[P]{P}{D. Popescu}, {\it Stanley depth of multigraded modules}, J. Algebra,
{\bf 321}(10), (2009),  2782-2797.

\bibitem[R]{R}{A. Rauf}, {\it Stanley decompositions, pretty clean filtrations and reductions
modulo regular elements},  Bull. Math. Soc. Sci. Math. Roumanie (N.S.), {\bf 50(98)}(4), (2007),
347-354.

\bibitem[S1]{S1} {A. Soleyman Jahan}, {\it Easy proofs of some well known facts via cleanness},
Bull. Math. Soc. Sci. Math. Roumanie, (N.S.), {\bf 54(102)}(3), (2011), 237-243.

\bibitem[S2]{S2}{A. Soleyman Jahan}, {\it Prime filtrations and Stanley decompositions
of squarefree modules and Alexander duality}, Manuscripta Math., {\bf 130}(4), (2009), 533-550.

\bibitem[S3]{S3}{A. Soleyman Jahan}, {\it Prime filtrations of monomial ideals and polarizations},
J. Algebra, {\bf  312}(2), (2007), 1011-1032.

\bibitem[SZ]{SZ}{A. Soleyman Jahan and X. Zheng}, {\it Monomial ideals of forest type},
Comm. Algebra, {\bf 40}(8), (2012), 2786-2797.

\bibitem[St]{St}{R.P. Stanley}, {\it Linear Diophantine equations and local cohomology}, Invent.
Math., {\bf 68}(2), (1982), 175-193.

\bibitem[TY]{TY}{N. Terai and K-I. Yoshida}, {\it Locally complete intersection Stanley-Reisner
ideals}, Illinois J. Math., {\bf 53}(2), (2009), 413-429.

\bibitem[V]{V}{R.H. Villarreal}, {\it Monomial Algebras}, Monographs and Textbooks in Pure and
Applied Mathematics, {\bf 238}, Marcel Dekker, Inc., New York, 2001.

\end{thebibliography}
\end{document}